\documentclass[aos,MSNbibl,nameyear,seceqn,dvips]{arximspdf}
\usepackage{dcolumn}

\doi{10.1214/12-AOS996} %kopijuoti is PTS
\volume{40}
\issue{3}
\pubyear{2012}
\firstpage{1465}
\lastpage{1488}

\makeatletter
\newtheorem{lemma}{Lemma}
\newtheorem{them}{Theorem}

\newcolumntype{d}[1]{D{.}{.}{#1}}
\makeatother

\begin{document}
\begin{frontmatter}

\title{Modeling left-truncated and right-censored survival data with
longitudinal covariates\thanksref{T1}}
\runtitle{JT modeling subject to LTRC}

\begin{aug}
\author[A]{\fnms{Yu-Ru} \snm{Su}\corref{}\ead[label=e1]{yrsu@stat.ncku.edu.tw}}
\and
\author[B]{\fnms{Jane-Ling} \snm{Wang}\ead[label=e2]{jlwang.ucdavis@gmail.com}}
\runauthor{Y.-R. Su and J.-L. Wang}
\affiliation{University of California, Davis and National Cheng Kung
University, and University of California, Davis}
\address[A]{Graduate Group in Biostatistics\\
University of California, Davis\\
Davis, California 95616\\
USA\\
and\\
Department of Statistics\\
National Cheng Kung University\\
Tainan 701\\
Taiwan\\
\printead{e1}} %adresu isvedimo komanda gale!
\address[B]{Department of Statistics\\
University of California, Davis\\
Davis, California 95616\\
USA\\
\printead{e2}}
\end{aug}
\thankstext{T1}{Supported in part by NIH Grant 1R01AG025218-01.}

% HISTORY:
\received{\smonth{3} \syear{2011}}
\revised{\smonth{3} \syear{2012}}

% ABSTRACT
%
\begin{abstract}
There is a surge in medical follow-up studies that include
longitudinal covariates in the modeling of survival data. So far,
the focus has been largely on right-censored survival data. We
consider survival data that are subject to both left truncation
and right censoring. Left truncation is well known to produce
biased sample. The sampling bias issue has been resolved in the
literature for the case which involves baseline or time-varying
covariates that are observable. The problem remains open, however,
for the important case where longitudinal covariates are present
in survival models.
%Traditional approaches, such as the partial likelihood approach for
%the Cox proportional hazards model, encounter difficulties when
%longitudinal covariates are involved in the modeling of the survival
%data.
A joint likelihood approach has been shown in the literature to
provide an effective way to overcome those difficulties for
right-censored data, but this approach faces substantial additional
challenges in the presence of left truncation. Here we thus
propose an alternative likelihood to overcome these difficulties
and show that the regression coefficient in the survival component
can be estimated unbiasedly and efficiently. Issues about the bias
for the longitudinal component are discussed. The new approach is
illustrated numerically through
simulations and data from a multi-center AIDS cohort study.
\end{abstract}

% KEYWORDS
%

\end{frontmatter}

%s1 ###
\section{Introduction}
Since the seminal paper by \citet{WulfT97}, longitudinal covariates
have played an increasingly important role in the modeling of
survival data. One major challenge to incorporate longitudinal
covariates is that simple approaches, such as the partial
likelihood method for the Cox proportional hazards model
[\citet{Cox72}], often require knowledge of the entire longitudinal
process. This is often not feasible in reality for follow-up
checks at discrete and intermittent time points. A common practice
is to impute the values of the missing longitudinal processes and
then apply the partial likelihood approach to the imputed data.
This is called a~two-stage approach, where the longitudinal
process is imputed at the first stage before the partial
likelihood approach is employed to estimate parameters in the
survival model at the second stage. The most common imputation
method is to use the last and most recent value of the patient to
impute a~missing value, the so-called last-value-carry-forward
method, which has been adopted in standard software such as SAS
and R. Additional two-stage procedures were developed by \citet
{TsiaDW95} and \citet{DafnT98}. %\citep{pawitan93,tsiatis95,dafni98}%.

It is easy to foresee serious biases with such an imputation
method if the follow-up schedule is infrequent over time and also
when the longitudinal covariates are contaminated by noises or
measurement errors. Both scenarios provide strong motivation to find
alternative approaches. The approach developed by \citet{WulfT97}
to model the survival and longitudinal data simultaneously through
their joint likelihood is attractive on two counts: (i) the
resulting parametric estimators are semiparametrically efficient
when the baseline hazard function is unknown, and (ii) the joint
likelihood procedure is often insensitive to the normality
assumption on the longitudinal data, if there is a reasonable
number of repeated measurements available for the longitudinal
processes; see \citet{ZengC05} and \citet{DupuGM06} for (i) and
\citet{SongDT02b}, \citet{TsiaD04} and \citet{HsieTW06} for (ii).

% [add those references and make sure the quotation is correct.]
% [ok!!]

The above joint likelihood approach not only successfully removes
the biases on the survival component but also leads to efficient
estimation. A~historical example for the joint likelihood approach
is the investigation of CD4 T-cell counts as a biomarker of
time-to-death or time-to-AIDS [\citet{DegrT94}, \citet{WulfT97}, \citet{HendDD00}]. In
these and other works, the survival time is subject to the usual
right censoring. However, left truncation is common for studies
with delayed entry. Specifically, if the recruitment of patients
continues after the onset time of a study, those that have already
experienced the event are often excluded from the study, which
then results in left truncation of the event-time. Patients who
remain in the study are further subject to the usual right
censoring, so the sample consists of left-truncated and right-censored
(LTRC) survival times. It is well known that left
truncation is a biased sampling plan as subjects with shorter
survival times tend to be excluded from the sample. As a result,
the longitudinal measurements are also sampled with bias.

An example of left-truncated and right-censored longitudinal study
is the Italian multi-center HIV (human immunodeficiency virus)
study [\citet{RezzLASPTSRAM89}, \citet{Ital92}], where the primary endpoint
is the time from HIV positive to AIDS onset, that is, the incubation
period of AIDS. In this study, patients who have developed AIDS at
the time of\vadjust{\goodbreak} recruitment were excluded from the study, resulting in
left truncation of the survival data, and CD4 counts for those who
were HIV positive but ADIS free were measured at each follow-up
visit. As there are no procedures available to handle such data
properly, we develop in this paper a~semiparametric joint likelihood approach
to accommodate LTRC survival data with longitudinal covariates
that are measured intermittently.

Although there is a sizable literature to jointly model
right-censored survival and longitudinal data
[see \citet{WulfT97}, \citet{HendDD00}, \citet{SongDT02b} and the review
papers by
\citet{TsiaD04}], %Davidian08%
the extension to LTRC survival data turns out nontrivial due to
the left-truncation feature of the data. To see this, consider
first the simpler case of left-truncated data with
time-independent covariates or no covariates at all.
\citet{Lynd71}, \citet{Wood85} and \citet{Wang87} investigated
estimation of the survival function when subjects come from the
same population, that is, there are no covariates involved. Here, one
only needs to adjust the risk set for truncated data to reach
a~suitable extension of the Kaplan--Meier estimator. For
time-independent covariates \citet{AndeBGK93} considered estimation
under the Cox model and showed that the partial likelihood
approach for right-censored data still works for LTRC survival
data when one conditions on the values of the covariates and
truncation times.

For time-dependent covariate, the Andersen et. al. (\citeyear{AndeBGK93}) partial
likelihood approach can still be employed if the entire covariate
history is available for all subjects. This is not the case for
longitudinal covariates that are observed intermittently at
discrete time points.
%
%Although imputation methods, such as the last-value-carry-forward
%approach in SAS and R packages, could be used to recover the
%entire history of the longitudinal process, the resulted
%prediction errors lead to bias in the estimates of the regression
%parameters much like the case of covariates measured with errors.
%
Since imputation methods lead to biases of the estimates, bias
corrected approaches have been employed in the literature for
right-censored data with longitudinal covariates. In particular,
\citet{Wang06} proposed a method to correct the bias through the
partial score equation. Such an approach is termed ``corrected
score'' methods, which originates from studies of measurement
errors. While corrected score methods typically lead to
$\sqrt{n}$-consistent estimators for the regression parameters in
the Cox model, they are not efficient and easy to derive.
Extension of the corrected score methods to LTRC (left-truncated
and right-censored) data might be feasible but have not been
explored. In this paper, we adopt the full and joint likelihood
approach of the survival and longitudinal data due to its
aforementioned efficiency and robustness features. Unfortunately,
direct maximization of the full joint likelihood is much more
complicated than the cases with no left truncation. We discovered
a modified likelihood that is simpler, yet retains the efficiency
of the full likelihood approach, as described in Section~\ref{sec2}.

The rest of the paper is organized as follows. In Section~\ref{sec2}, we
introduce a~joint model setting for both the survival time and
longitudinal processes and propose a modified likelihood approach
for statistical inference. An EM algorithm to maximize the
modified likelihood is derived in Section~\ref{sec3}, along with the large\vadjust{\goodbreak}
sample properties of the nonparametric maximum modified likelihood
estimator (NPMMLE), including consistency, asymptotic normality
and efficiency. Numerical performance of the proposed estimating
procedure is validated through simulation studies in Section~\ref{sec4} and
illustrated through the Italian HIV study in Section~\ref{sec5}.
Section~\ref{sec6}
contains some discussion.\vspace*{-2pt}

%s2 ###
\section{Joint modeling under LTRC}\label{sec2}

We consider the setting that the survival time $Y^*$ of a subject
is subject to random left truncation by $T^*$, so a~subject is
enrolled in a study only if $Y^*\geq T^*$. Let $n$ be the total
number of subjects enrolled in the study. With such a biased
sampling plan, to avoid confusion of notation, we denote the
survival and truncation time of the $i$th enrolled subjects as
$(Y_i, T_i)$, which are sampled from the joint subpopulation of
$(Y_i^*, T_i^*)$, where $Y_i^* \geq T_i^*$. Upon entering the
study, these $n$ subjects are subject to the usual right
censorship, so the final observed survival data for the $i$th
subject is a triplet $(T_i, Z_i, \Delta_i)$, where
$Z_i=\min(Y_i,C_i)$ is the time of the endpoint event or drop-out
(censoring) time $C_i$, whichever occurs first, and $\Delta_i=I
(Y_i \leq C_i)$ is the censoring indicator.

In reality, drop-out or censoring only occurs when a subject is
enrolled into
the study. This fact implies that the right-censoring time $C_i$
is greater than the truncation time $T_i$, for $i=1,\ldots,n$.
Therefore, we introduce a positive random variable $U_i$ to
represent the time from entry into the study to drop-out from the
study, that is, $U_i=C_i-T_i$.
%Adopting this decomposition of the censoring time would benefit us to
%derive the likelihood and get rid of the distribution of the censoring
%time under some proper assumptions which will be demonstrated later.

In addition to the survival data, baseline and longitudinal
covariates are collected intermittently for the $i$th subject from
the time the subject enters the study until the observational
limit $Z_i$. This results in $n_i$ repeated measurements, denoted
by $\vec{W}_i=(W_{i1},W_{i2},\ldots,W_{in_i})$, where the
measurements are taken at time points
$\vec{s}_i=(s_{i1},s_{i2},\ldots,s_{in_i})$. It is important to
make a~note here that the observed $\vec{W}_i$ are also subject to
the same biased sampling plan as the survival data, so there is a
background longitudinal vector, which we will denote as
$\vec{W}^*_i$ for the $i$th subject enrolled in the study.
Therefore,~$\vec{W}_i$ is sampled from the subpopulation of
$\vec{W}^*$,  where $Y_i^* \geq T_i^*$, and values beyond~$Z_i$ are
not observed. For simplicity of notation, we assume in this
section that there is only one longitudinal covariates, but
additional longitudinal or baseline covariates can be handled
easily and the AIDS data discussed in Section~\ref{sec5} contain two
longitudinal covariates, one observed intermittently but the
complete history of the other one, the time-dependent treatment
indicator, is available.

%
%the focus lies on the association between an event-time $Y^*_i$,
%and a longitudinal process $X_i(\cdot)$. Since some subjects may
%be lost to follow up during the time course or never encounter the
%event of interest by the end of the study, what one can truly
%observed are and $\Delta_i=I(Y^*_i\leq C_i)$, where $C_i$ stands
%for the right-censoring time. Besides the issue of
%right-censoring, the data is also subjected to the problem of
%left-truncation with $T_i$ as the truncating time. Thus the
%observed survival data from the $i$th subject is a triple as
%$(T_i,Y_i,\Delta_i)$ only when $Y^*_i>T_i$.
%
%Meanwhile, a longitudinal process is measured intermittently $n_i$
%
%times since the $i$th subject entered the study till the
%occurrence of either the drop-out or the survival event.
%

%s2.1 ###
\subsection{The joint models}

Since repeated measurements from the same subjects are likely to
be correlated, we introduce a latent $q\times1$ random
vector~$A^*_i$ to account for their dependency and assume a common
parametric density function $f_A^*(\cdot\vert\alpha)$ with an
unknown parameter $\alpha$ for $A_i^*$. A linear mixed effects
model will be considered for the longitudinal covariate
%
%e2.1 ###
\begin{equation}
\vec{W}^*_i=X(\vec{s}_i)+\varepsilon_i= g(\vec{s}_i)A^*_i+\varepsilon_i,\vadjust{\goodbreak}
\end{equation}
where $g(\cdot)$ is a known $q$-dimensional function and the
$n_i\times1$ vector $\varepsilon_i$ plays the role of measurement
errors, sampled from a multivariate normal distribution
with independent marginal distribution $\mathcal{N}(0,\sigma^2)$,
and independent of all other aforementioned random variables.

% will assume that the survival time is associated with the
%longitudinal covariate only through the latent variable $A_i$.

For the survival time $Y^*_i$, a proportional hazards model is
employed, and the hazard rate of $Y_i^*$ at time $t$ given $A^*_i$
is
%
%e2.2 ###
\begin{equation}
\lambda_{Y^*_i}(t\vert A^*_i)=\lambda_0(t)\exp(\beta X_i(t)),
\end{equation}
where $\lambda_0$ is the baseline hazard rate and $\beta$ is the
regression coefficient. The truncation time $T^*_i$ and the time $U_i$,
from entry to drop-out, are assumed to have distribution
function $F_{T^*}(\cdot)$ and $F_U(\cdot)$, respectively. We adopt
the standard assumption in survival analysis, that $Y^*_i$,
$T^*_i$ and $U_i$ are conditionally independent given the
covariates. This is equivalent to assuming conditional
independence of $Y^*_i$, $T^*_i$ and $U_i$ given the value of
$A^*_i$. We also assume that $T^*_i$ and $U_i$ are independent of
$A^*_i$, and
%[Jane-Ling -- I think the sentence below should be $T^*_i$ and
%$U_i$ are noninformative for the parameters in the models.....]
the parameters in the models for either the survival or
longitudinal parts are noninformative.\vspace*{-2pt}

%
%We also assume that the truncation and censoring scheme are
%noninformative, given the value of $A_i$.
%

%s2.2 ###
\subsection{A modified likelihood approach}
% Observed likelihood%

%[The notations $t_i$ for observed truncation time and $\vec{t}_i$ for
%the longitudinal schedule is confusing, so %I'm considering to replace
%all $\vec{t}_i$ by $\vec{s}_i$. If you agree. make the switch
%throughout the paper. The switch has been made.]

For the model described in the previous subsection, the parameters
of interest are ($\beta$, $\alpha$, $\sigma^2$ and
$\Lambda_0(\cdot)$), where the first three components are in the
Euclidean space whereas $\Lambda_0(t)= \int_0^{t} \lambda_0(u) \,
du$, the cumulative hazard function, is in a functional space,
hence the model is semiparametric. Since a likelihood approach
usually provides the most efficient estimating procedure, we first
consider the full likelihood function $L_i^O$ based on the
observations ($t_i$, $z_i$, $\delta_i$, $\vec{w}_i$) from the
$i$th subject. The derivation of the full likelihood from the
$i$th subject is shown below:
%
%e2.3 ###
\begin{eqnarray}\label{originalfulllike}
L_i^O &= & f_{(T,Y,\Delta,W)}(t_i,z_i,\delta_i,w_i)
= \frac{f_{(T^*,Y^*,\Delta^*,W^*)}(t_i,z_i,\delta_i,w_i)}{P(Y^*\geq
T^*)}\nonumber
\\[-1pt]
&=&\biggl\{\int[f_{Y^*}(z_i\vert A^*_i=a_i)]^{\delta_i}
[S_{Y^*}(z_i\vert A^*_i=a_i)]^{1-\delta_i}\nonumber
\\[-1pt]
&&\hspace*{55pt}{}\times
f_{W^*}(w_i\vert A^*_i=a_i) f_{A^*}(a_i)\, da_i\biggr\}
f_{T^*}(t_i)/{\mathrm{P}( Y^*_i \geq T^*_i)}\nonumber
\\[-1pt]
&=&\biggl\{ \int\frac{[f_{Y^*}(z_i\vert A^*_i=a_i)]^{\delta_i}
[S_{Y^*}(z_i\vert A^*_i=a_i)]^{1-\delta_i}}{S_{Y^*}(t_i|A^*_i=a_i)}\nonumber
\\[-8pt]
\\[-8pt]
&&\hphantom{\biggl\{ \int}{}\times f_{W^*}(w_i\vert A^*_i=a_i)\frac{S_{Y^*}(t_i\vert A^*=a_i)f_{A^*}(a_i)}{S_{Y^*}(t_i)}
\,da_i\biggr\}
\frac{S_{Y^*}(t_i)f_{T^*}(t_i)}{\mathrm{P}( Y^*_i \geq
T^*_i)}\nonumber
\\[-1pt]
&=& \biggl\{\int[f_{Y^*}(z_i\vert Y_i^*\geq t_i, A^*_i=a_i)]^{\delta_i}
[S_{Y^*}(z_i\vert Y^*_i\geq t_i, A^*_i=a_i)]^{1-\delta_i}\nonumber
\\[-1pt]
&&\hspace*{95pt}{}\times f_{W^*}(w_i\vert A^*_i=a_i) f_{A^*}(a_i\vert Y^*_i\geq t_i)\, da_i\biggr\}\nonumber
\\[-1pt]
&&{}\times f_{T^*}(t_i\vert Y^*_i \geq T^*_i),\nonumber
\end{eqnarray}
where $f_V$ is the density function of the random variable $V$ in
the subscript, and $S_V$ is the corresponding survival function.
In (\ref{originalfulllike}), besides the baseline hazard function
$\lambda_0$, the density function $f_{T^*}$ also serves as a
nonparametric component. Because of these two nonparametric
components, the full likelihood function is unbounded, so we
resort to the nonparametric maximum likelihood approach, which
leads to a similar scenario as in conventional survival analysis
that the full likelihood is the same as the conditional likelihood
given the left-truncation time. This has been explored in the
literature [\citet{AndeBGK93}, \citet{KleiM03}] for LTRC data with baseline
covariates, and was first explored in \citet{Wang87} for
the simpler situation of left-truncated data that came from a
single population. Following a similar argument as in
\citet{Wang87}, we found that the full likelihood can be simplified
to the following conditional likelihood for the $i$th subject~as
%
%e2.4 ###
\begin{eqnarray}\label{fulllike}\qquad
L_i^C
&=& \biggl\{\int[f_{Y^*}(z_i\vert Y_i^*\geq t_i, A^*_i=a_i)]^{\delta_i}
[S_{Y^*}(z_i\vert Y^*_i\geq t_i, A^*_i=a_i)]^{1-\delta_i}\nonumber
\\[-10pt]
\\[-10pt]
&&\hspace*{93pt}{}\times f_{W^*}(w_i\vert A^*_i=a_i) f_{A^*}(a_i\vert Y^*_i\geq t_i) \,
da_i\biggr\}.\nonumber
\end{eqnarray}
Next, we consider the nonparametric maximum likelihood estimators\break
(NPMLE) of the survival component, which, by a similar argument
for joint modeling right-censored data and their longitudinal
covariates [\citet{ZengC05}, \citet{DupuGM06}], leads to a piecewise linear
baseline cumulative hazard function with jumps at each uncensored
event time (i.e., at $Y_i$, whenever $\Delta_i=1$). Let $n_u$
denote the total number of uncensored events, the baseline
cumulative hazard function is thus re-parameterized as a
$n_u$-dimensional vector.

%
%the likelihood function is unbounded when the baseline hazard
%function is unrestricted, or when the truncation distribution is
%modeled nonparametrically. Therefore, we consider instead the
%nonparametric maximum likelihood estimators (NPMLE), which
%correspond to a piecewise linear baseline cumulative hazard
%function with jumps at each uncensored event time (i.e. at $Y_i$,
%whenever $\Delta_i=1$). Let $n_u$ denote the total number of
%uncensored events, the baseline cumulative hazard function is thus
%re-parameterized as a $n_u$-dimensional vector.
%
%%\iffalse%%%%%%%%%%%%%
%
%Tracing back to the conventional Cox model when the entire history
%of time-dependent covariates is observed, the partial likelihood
%for right censored data can be extended to LTRC data by
%considering the conditional likelihood given the left-truncation
%time and the covariates, an approach that is common in the
%literature \citep{AndeBGK93,KleiM03}. It turns out that the NPMLE
%of the survival-related parameters, $\beta$ and $\lambda_0$, based
%on the full likelihood, is the same as the conditional likelihood
%given the left-truncation time and the covariates. This can be
%seen from a similar argument as in \citet{Wang87}, who dealt with
%the simpler situation for left truncated data that came from a
%single population. Such a derivation of the NPMLE based on full
%likelihood is also similar to the much investigated case of a
%joint modeling setting with right censored data, where NPMLE for
%the parametric components enjoy nice asymptotic properties and are
%semiparametrically efficient.
%

So far, the derivation of the likelihood function and NPMLE
follows a similar path as the much investigated case of a joint
modeling setting with right-censored data, where NPMLE's for the
parametric component enjoy nice asymptotic properties and are
semiparametrically efficient. Despite these similarities, the left
truncation feature triggers complications in the estimation of the
finite dimensional parameter in the joint LTRC model. First, as
shown in the Appendix, the parameter $\alpha$ associated with the
latent variable $A^*$ is not identifiable. This is a consequence
of the biased sampling plan, since the samples are actually drawn
from the subpopulation $Y^* \geq T^*$. Consequently, only
$E(A^*|Y^*\geq T^*)$ and $\mathrm{var}(A^*|Y^* \geq T^*)$ could be
identified under the normality assumption. Thus, while it is
possible to identify the unknown parameters of $Y^*$ and $T^*$
based on the joint conditional distribution of $(Y^*,T^*)| Y^*
\geq T^*$, where the notation $(\cdot|Y^*\geq T^*)$ stands for a
random variable/vector sampled from the subpopulation with
$Y^*\geq T^*$, there is not enough information to recover $E(A^*)$
and var$(A^*)$ and hence the true longitudinal parameters
$\alpha$.

A second complication is that the score equations for the survival
components, $\beta$ and $\Lambda_0$,\vadjust{\goodbreak}
%$(\lambda_1,\ldots,\lambda_{n_u})$.
are much more complicated than the situation under a~right-censored
only model and, as shown in Appendix~\ref{appa1}, as they require
estimation of the expectations of nonlinear functions of the
observed data along with the the parameters of interest. This
motivates us to modify the likelihood so as to simplify the
estimation of all parameters that are identifiable. Our proposal
is to aim at the following \textit{modified likelihood}, denoted
by $L^{m}$, as an alternative of the full, also the conditional,
likelihood in (\ref{fulllike}). The modified likelihood is
%
%e2.5 ###
\begin{eqnarray}\label{marparlike}
L^{m} &=& \prod_{i=1}^{n}\biggl\{ \int[f_{Y^*}(z_i\vert
Y^*_i\geq t_i, A^*_i=a_i)]^{\delta_i}\nonumber
\\
&&\hphantom{\prod_{i=1}^{n}\biggl\{ \int}{}\times[S_{Y^*}(z_i\vert Y^*_i\geq t_i, A^*_i=a_i)]^{1-\delta_i}
\\
&&\hspace*{34pt}{}\times f_{W^*}(w_i\vert A^*_i=a_i) f_{A^*}(a_i)\,
da_i\biggr\}f_{T^*}(t_i|Y^*_i\geq T^*_i),\nonumber
\end{eqnarray}
where the lower case variables denote the values of the
corresponding upper case variables, for example, $\delta_i$ is the value
of $\Delta_i$. The estimators obtained by maximizing the modified
likelihood, where the nonparametric cumulative hard function is
replaced by a
step function will be referred to as the \textit{nonparametric maximum
modified likelihood} (NPMMLE) hereafter.

%Similar to the idea of NPMLE, we replace the baseline hazard function $
%representing the jump sizes of a piecewise constant cumulative hazard
%function.

The difference between (\ref{fulllike}) and (\ref{marparlike}) is
that $f_{A^*}(a_i \vert Y^*_i\geq t_i)$ in the full likelihood
(\ref{fulllike}) is replaced by $f_{A^*}(a_i)$ in
(\ref{marparlike}). This is motivated by the fact that
$f_{A^*}(a\vert Y^*\geq t)=\frac{S_{Y^*}(t\vert
A^*=a)}{S_{Y^*}(t)}f_{A^*}(a)$ and $E[\frac{S_{Y^*}(t\vert
A^*)}{S_{Y^*}(t)}]=1$, for any
% ==== Lemma 1 might be moved to the content instead of being in the
%appendix ==== %
$t$, and that, as shown in Lemma~\ref{lem1} in the \hyperref[app]{Appendix}, the score
%
%================================================================================
%%
functions of the survival parameters from (\ref{marparlike}) are
asymptotically the same as those from (\ref{fulllike}).
%This suggests the validity of this modified likelihood.
Theoretical results in the next section and numerical evidence in
Section~\ref{sec4} demonstrate good performance of estimators of all the
survival parameters, $(\beta,\Lambda_0 (\cdot))$ and of the
measurement errors $\sigma^2$ of the longitudinal component that
we derived from this modified likelihood.

\section{EM-algorithm and asymptotic properties}\label{sec3}

Let $\gamma= (\beta, \alpha, \sigma^2)$ be the finite dimensional
parameter in the joint survival and longitudinal model, and
$\Lambda$ be a step function. The log modified likelihood is
%
%[Jane-Ling -- Should we change the notations below for the sorted
%observed survival times from $y_j^0$ to $z_j^0$?
%
%no, we still keep $y_j^0$ since they are uncensored ones for which
%$y$ is the same as $z$.]
%
\begin{eqnarray*}
l^{m}(\gamma,\Lambda)&=&\sum_{i=1}^{n} \ln
\int[\Lambda\{z_i\}\exp{\beta g(z_i)a_i}]^{\delta_i}
\\
&&\hphantom{\sum_{i=1}^{n} \ln
\int}{}\times \exp\biggl\{-\sum_{j:t_i<y_j^0\leq z_i}\Lambda\{y_j^0\}\exp\{\beta g(y_j^0)a_i\}\biggr\}
 \\
&&\hphantom{\sum_{i=1}^{n} \ln
\int} {}\times(2\pi
\sigma^2)^{-n_i/2}
\\
&&\hphantom{\sum_{i=1}^{n} \ln
\int}{}\times\exp\Biggl\{-\sum_{j=1}^{m_i}
[w_{ij}-g(s_{ij})a_i]^2/(2\sigma^2)\Biggr\}
f_{A^*}(a_i) \, da_i,\nonumber
\end{eqnarray*}
where $\Lambda\{\cdot\}$ is the jump size of $\Lambda$ at the
respective time point in the argument, and $y_j^0$ is the $j$th
sorted observed survival time in increasing order.
%We will assume that a total of $n_u$ such survival times are observed.
Moreover,~$\tau_1$ and $\tau_2$ denote the lower bound of
truncation time and the largest censoring time corresponding to
the end of the study.

Since direct maximizing the proposed modified likelihood involves
integration of a complex function with respect to the random
effects, we employ the expectation-maximization (EM) algorithm
[\citet{LairW82}] to stabilize the maximization procedure. In the
implementation of the EM algorithm, a Monte Carlo integration
approach is used to approximate the expectation terms of functions
$h(A^*)$ appearing in the E-step.
%according to the posterior distribution given the observed data.
A one-step Newton--Raphson method is applied to solve the nonlinear
equations in the M-step. The posterior density of the random
effects $A^*_i$ given the observed data from the $i$th subject,
$o_{i}= (t_i, z_i, \delta_i, \vec{w}_i)$, is of the form
%
%%%%%%%%%%%%%%%%%%%%%%%%%%%%%%%
%[ Shouldn't all the $w_i$ be $\vec{w}_i$ and $w$ be $\vec{w}$?
%Correct this throughout the paper. You have not defined $y_i^0$
%even though it is clear what it should be. This kind of oversights
%happened repeatedly throughout the paper. I've spotted and
%corrected some but there are probably more that needs to be
%caught. Another example is the $\hat{\Lambda}_k$ below. While it
%is clear what the other hats are, $\Lambda_k$ has never been
%defined. Pls confirm that this is for the cumulative hazard
%function, not the hazard function. Please do a THOROUGH check of
%the entire paper. This is a common mistake for first time writer,
%and it gets better with time - but do check carefully in the
%future.
%
%Thanks for pointing this out. I shall be more careful.]
%
%%%%%%%%%%%%%%%%%%%%%%%%%%%%%%%%%
%
% ==== Note: 2011.2.15 ===== %
% The validity of the posterior density is based on the proof in Fisher
%consistency. Please see the sketch.
% The notation A^* below should be A which is A^*|Y^*>T^* because the
%observed data is from the truncated population.
% To be consistent with the A^* in the last two sentence above, we keep
%using the notation A^*.
% The first equation is based on the form in the second equation.
% Actually, A^*|W^* is the same as A|W if we impose an assumption:
%f_{A}=f_{A^*}.
% This assumption makes the full likelihood become modified likelihood.
%
\begin{eqnarray*}
f_{A^*|O}(a|o_i)&=&\frac{f_{(Y,\Delta)|(A,T)}(z_i,\delta_i|a,t_i)\times
f_{A^*|W^*}(a|\vec{w}_i)}{\int
f_{(Y,\Delta)|(A,T)}(z_i,\delta_i|a,t_i)\times
f_{A^*|W^*}(a|\vec{w}_i)\,da}
\\
&=&[\Lambda\{z_i\}]^{\delta_i}\exp\biggl\{-\sum_{j:t_i<y_j^0\leq
z_i} \Lambda\{y_j^0\}\exp\{\beta g(y_j^0)a\}\biggr\}\times
f_{A^*|W^*}(a|\vec{w}_i)
\\
&&{}\Big/\int
[\Lambda\{z_i\}]^{\delta_i}\exp\biggl\{-\sum_{j:t_i<y_j^0\leq
z_i}\Lambda\{y_j^0\}\exp\{\beta g(y_j^0)a\}\biggr\}
\\
&&\hphantom{{}/\int}{}\times
f_{A^*|W^*}(a|\vec{w}_i)\,da.
\end{eqnarray*}
For a simpler implementation of the algorithm, we shall impose a
normal assumption on the random effects and assume that $A^*_i$,
$i=1,\ldots,n,$ follow a~normal distribution $N(\mu,\Sigma)$,
where $(\mu,\Sigma)$ plays the role of the parameter~$\alpha$.

By taking the first derivative of the log modified likelihood
calculated in the E-step\vspace*{1pt} with respect to each parameter, the
NPMMLE, $\hat{\beta}$, $\{\hat{\lambda}_k,k=1,\ldots,n_u\}$,
$\hat{\sigma}$, $\hat{\mu}$ and $\hat{\Sigma}$, can be obtained
through the following formulas, where $\lambda_k$ is the jump size
of $\Lambda$ at the $k$th sorted observed survival time:
\begin{eqnarray*}
\hat{\lambda}_k&=&\frac{1}{\sum_{i:t_i<y_k^0\leq
z_i}\mathrm{E}[\exp\{\beta
g(y_k^0)A^*_i\}|o_i]},\qquad k=1,\ldots,n_u,
\\
\hat{\sigma}^2&=&\frac{1}{\sum_{i=1}^{n} n_i}\sum_{i=1}^{n}\sum_{j=1}^{n_i}
E\bigl[\bigl(w_{ij}-g(s_{ij})A^*_i\bigr)^2|o_i\bigr],
\\
\hat{\mu}&=&\frac{1}{n}\sum_{i=1}^{n}
E(A^*_i|o_i),\qquad\hat{\Sigma}=\frac{1}{n}\sum_{i=1}^{n}E[(A^*_i-\hat
{\mu})(A^*_i-\hat{\mu})^T|o_i]
\end{eqnarray*}
and $\hat{\beta}$ is the root of the score $s(\beta)$, which is
solved by an one-step Newton--Raphson method with the updating rule
\[
\beta_{\mathrm{new}}=\beta_{\mathrm{old}}-\frac{s(\beta_{\mathrm
{old}})}{s'(\beta_{\mathrm{old}})},
\]
where
\begin{eqnarray*}
s(\beta)&=&\sum_{i=1}^{n}\delta
_i\biggl[g(z_i)E(A^*_i|o_i)-\frac{\sum_{j:t_j<z_i\leq
z_j}E(g(z_i)A^*_j\exp\{\beta
g(z_i)A^*_j\}|o_j)}{\sum_{j:t_j<z_i\leq
z_j}E(\exp\{\beta g(z_i)A^*_j\}|o_j)}\biggr],
\\
s'(\beta)&=&\sum_{i=1}^{n}\delta_i\biggl\{\biggl[\frac{\sum_{j:t_j<z_i\leq
z_j}E(g(z_i)A^*_j\exp\{\beta
g(z_i)A^*_j\}|o_j)}{\sum_{j:t_j<z_i\leq
z_j}E(\exp\{\beta g(z_i)A^*_j\}|o_j)}\biggr]^2
\\
&&\hphantom{\sum_{i=1}^{n}\delta_i\biggl\{}{}-\frac{\sum_{j:t_j<z_i\leq
z_j}E((g(z_i)A^*_j)^2\exp\{\beta
g(z_i)A^*_j\}|o_j)}{\sum_{j:t_j<z_i\leq
z_j}E(\exp\{\beta g(z_i)A^*_j\}|o_j)}\biggr\}.
\end{eqnarray*}

Except for $\alpha$, the proposed nonparametric maximum modified
likelihood estimates (NPMMLE) of the parameters enjoy nice
properties that are similar to the NPMLE, as illustrated in the
next two theorems. Below, we list some regularity conditions
needed for the theorems:
\begin{enumerate}[(C6)]
\item[(C1)] The parameter space of the finite dimensional parameters,
$S_{\gamma}$, is bound\-ed and closed on Euclidean space. The
true value $\gamma_0$ is an interior point of $S_{\gamma}$.
\item[(C2)] On the parameter space of $\beta$, $(\exp\{\beta g(S)A^*\}
|Y^*\geq T^*)$ is bounded below by $m$ and above by $M$ with
probability 1.
\item[(C3)] $P(T\leq\tau_1\mbox{ and }Y\geq\tau_2 )> 0$.
This ensures that not all data are truncated or censored.
\item[(C4)] $E_{\theta_0}\{\exp[\beta_0g(u)A^*]I(T^*<u\leq Y^*)|Y^*\geq
T^*\}$ is bounded away
from 0 on the parameter space of $\beta$. Here
$E_{\theta_0}(\cdot)$ stands for the expectation taken under
the true value of the parameter $\theta_0$.
\item[(C5)] $g(t)$ is of uniformly bounded variation on
$[\tau_1,\tau_2]$, and there exists a~constant $D$ such that
$P(n_i\leq D)=1, \forall i$.
%$n_i$ the total number of measurements for each
%subject is bounded above by a certain $D$ with probability 1.
%
\item[(C6)] The distribution $f_{A^*}(\cdot|\alpha)$ is continuous with
respect to $\alpha$ and has continuous second derivative with
respect to $\alpha$. Moreover, the Fisher information matrix
obtained from $f_{A^*}$ for $\alpha$ is positive definite.
\end{enumerate}

%We illustrate this point with the following two asymptotic properties
%of the NPMMLE.

\begin{them}[(Consistency of the estimators)]
Under the regularity conditions \textup{C1--C5}, the NPMMLE of
$(\beta_0,\sigma^2_0,\Lambda_0)$, denoted as
$(\hat{\beta}_n,\hat{\sigma}^2,\hat{\Lambda}_n)$, is consistent
under the Euclidean norm $|\cdot|$ and supremum norm $\Vert
\cdot\Vert_{\infty}$ on $[\tau_1,\tau_2]$, respectively.
\end{them}

For $H=\{h=(h_1,h_2,h_3)\}$ and $0<p<\infty$, let $H_p=\{h\in
H\dvtx\Vert h_1\Vert, \vert h_2\vert,\allowbreak \Vert
h_3\Vert_v\leq p\}$, be a collection of directions that are
used in the \hyperref[app]{Appendix}. The notation $\Vert\cdot\Vert_v$
denotes the the total variation
of the function in the norm plus the absolute value of this function
evaluated at 0. The next theorem shows that the NPMMLE converges
in distribution to a Gaussian element in the parameter space at a
$\sqrt{n}$-rate.

\begin{them}[(Asymptotic normality and efficiency)]
Under the regularity conditions \textup{C1--C6}, the process
$\sqrt{n} (\hat{\alpha}_n-\mathrm{E}(\hat{\alpha}_n),\hat{\sigma
}^2_n-\sigma^2_0,\hat{\beta}_n-\beta_0,\hat{\Lambda}_n-\Lambda_0 )$
converges in distribution to a mean zero Gaussian process $G$ in the
functional space $l_{\infty}(H_p)$ on
$H_p$. Moreover, the NPMMLE $\hat{\beta}$ is semiparametrically
efficient for $\beta_0$.
\end{them}

Proofs of these two theorems are provided in the
\hyperref[app]{Appendix}.

For estimating the standard errors of the NPMMLE, we recommend to
use the bootstrap procedure instead of the profile likelihood
approach in \citet{MurpV00} and \citet{ZengC05}, which did not work
well for LTRC data due to the high fluctuation of the estimated
profile likelihood function and possibly negative estimate of the
standard error. The performance of the bootstrap procedure for
estimating the standard errors of the NPMLE under joint modeling
with right-censoring cases has been studied by \citet{TsenHW05}
for the accelerated failure time model, and by \citet{HsieTW06} for
the Cox model. The results in these two papers and support the validity
of the bootstrap method in the scope of joint modeling. Our simulation
results reported in Section~\ref{sec4} also supports the use of the bootstrap approach.
In
comparison, the bootstrap method is more reliable than the profile
likelihood method at a higher
computational cost.

%s4 ###
\section{Simulation study}\label{sec4}
To verify numerically the validity of the proposed procedure, we
conducted simulations under five
different settings. Since there is an intrinsic bias on the
longitudinal component, the simulations focus on the performance
of the estimate of $\beta$ and how it would be affected by the
level of contamination from the measurement errors and the
variation of the random effects. As a benchmark setting, we
considered a linear trend in time with random effects on the
longitudinal covariate and assess the influence of the variance of
the random slope on the accuracy of estimating~$\beta$. The
left-truncation times are generated from an exponential
distribution with parameter 1, while the right-censoring times are
from an exponential distribution with parameter 3. The baseline
hazard rate is from an exponential distribution with mean 1. All 5
simulation settings have sample size $n=200$ with true values
$\beta=1$, $\mu= (2,0.5)$ and
$(\sigma_{11},\sigma_{12})=(0.5,-0.001)$. The values of
$(\sigma_{22},\sigma^2)$ are different for the five settings and
set as: $(0.01,0.1)$, $(0.01,0.4)$, $(0.01,0.025)$, $(0.0025,0.1)$ and
$(0.04,0.1)$. The first three settings demonstrate the impact of
contaminations by measurement errors while the last two illustrate
the effect of the variation of the random slope.\vadjust{\goodbreak}

Simulation results based on 100 Monte Carlo samples are reported
in Table~\ref{tab1}. Results under the first three settings suggest that
$\beta$ can be estimated unbiasedly, and measurement errors affect
the precision, but not the magnitude of the biases. As expected,
higher level of noise contamination leads to less precise estimate
of $\beta$ and higher chance of divergence in the algorithm. In
all three settings, the variance of measurement errors can be
estimated with high accuracy and precision. Comparing with the
results under the first, fourth and fifth setting from Table~\ref{tab1}, we
observe that the variance of the random slopes has little effect
on the performance of $\hat{\beta}$.

%t1 ###
\begin{table}[t!]
\tablewidth=\textwidth
\tabcolsep=0pt
\caption{Simulation results under five settings with
sample size 200 and varying values of~$\sigma_{22}$
and~$\sigma^2$. The actual targets of the longitudinal estimates are
conditional quantities marked as~$\mu_1^*$ and~$\mu_2^*$ etc.
and are listed next to the
true longitudinal value in the first column. The mean and SD of the
estimates based on 100 Monte Carlo samples are reported in the second
and third column}\label{tab1}
\begin{tabular*}{\textwidth}{@{\extracolsep{\fill}}lcd{2.4}cd{2.4}c@{}}
\hline
&&\multicolumn{1}{c}{\textbf{Average of}}\\
\textbf{Case} & \textbf{Parameter} & \multicolumn{1}{c}{\textbf{NPMMLE}} & \multicolumn{1}{c}{\textbf{SE(MC)}}
& \multicolumn{1}{c}{\textbf{MSE}} & \textbf{Convergence rate}\\
\hline
1 & $\beta$ (1) & 0.9923 & 0.1633\phantom{$\mathrm{e}\mbox{-}4$} & 0.0267 & 98\%\\
% PL: 0.1430 (53/100)
& $\sigma^2$ (0.1) & 0.0998 & 0.0021\phantom{$\mathrm{e}\mbox{-}4$} & \multicolumn{1}{c}{$5\mathrm{e}{-}6$}&\\
& $\mu_1/\mu^*_1$ (2$/$1.73) & 1.7461 & 0.0478\phantom{$\mathrm{e}\mbox{-}4$} & 0.0668 &\\
& $\mu_2/\mu^*_2$ (0.50$/$0.50) & 0.4545 & 0.0985\phantom{$\mathrm{e}\mbox{-}4$} & 0.0118 &\\
& $\sigma_{11}/\sigma^*_{11}$ (0.50$/$0.45) & 0.4634 & 0.0527\phantom{$\mathrm{e}\mbox{-}4$} & 0.0041 &\\
& $\sigma_{12}/\sigma^*_{12}$ ($-$0.001$/$$-$0.001) & -0.0424 & 0.0453\phantom{$\mathrm{e}\mbox{-}4$} & 0.0038 &\\
& $\sigma_{22}/\sigma^*_{22}$ (0.01$/$0.01) & 0.0738 & 0.0409\phantom{$\mathrm{e}\mbox{-}4$} & 0.0057
&\\[5pt]
% Mean and Cov from the truly simulated random effects:
% (1.7321,0.5016) & (0.4547,-0.0009,0.0100)
2 & $\beta$ (1) & 0.9185 & 0.1765\phantom{$\mathrm{e}\mbox{-}4$} & 0.0378 & 72\%\\
% PL: 0.2309 (46/100)
& $\sigma^2$ (0.4) & 0.4003 & 0.0086\phantom{$\mathrm{e}\mbox{-}4$} & \multicolumn{1}{c}{$7\mathrm{e}{-}5$} & \\
& $\mu_1/\mu^*_1$ (2$/$1.74) & 1.7455 & 0.0531\phantom{$\mathrm{e}\mbox{-}4$} & 0.0676 & \\
& $\mu_2/\mu^*_2$ (0.50$/$0.50) & 0.3801 & 0.1640\phantom{$\mathrm{e}\mbox{-}4$} & 0.413 & \\
& $\sigma_{11}/\sigma^*_{11}$ (0.5$/$0.45) & 0.4730 & 0.0505\phantom{$\mathrm{e}\mbox{-}4$} & 0.0033 & \\
& $\sigma_{12}/\sigma^*_{12}$ ($-$0.001$/$$-$0.001) & -0.1122 & 0.0917\phantom{$\mathrm{e}\mbox{-}4$} &
0.0208 & \\
& $\sigma_{22}/\sigma^*_{22}$ (0.01$/$0.01) & 0.1856 & 0.1156\phantom{$\mathrm{e}\mbox{-}4$} & 0.0442 &
\\[5pt]
% Mean and Cov from the truly simulated random effects:
% (1.7359,0.5002) & (0.4497,-0.0011,0.0099)
3 & $\beta$ (1) & 1.0380 & 0.1548\phantom{$\mathrm{e}\mbox{-}4$} & 0.0254 & 96\%\\
% PL: 0.1276 (49/100)
& $\sigma^2$ (0.025) & 0.0250 & $4.8283\mathrm{e}\mbox{-}4$ & \simeq 0& \\
& $\mu_1/\mu^*_1$ (2$/$1.73) & 1.7443 & 0.0468\phantom{$\mathrm{e}\mbox{-}4$} & 0.0676 & \\
& $\mu_2/\mu^*_2$ (0.50$/$0.50) & 0.4900 & 0.0643\phantom{$\mathrm{e}\mbox{-}4$} & 0.0042 & \\
& $\sigma_{11}/\sigma^*_{11}$ (0.50$/$0.45) & 0.4520 & 0.0534\phantom{$\mathrm{e}\mbox{-}4$} & 0.0052 &
\\
& $\sigma_{12}/\sigma^*_{12}$ ($-$0.0004$/$$-$0.0004) & -0.0193 & 0.0338\phantom{$\mathrm{e}\mbox{-}4$} &
0.0015 & \\
& $\sigma_{22}/\sigma^*_{22}$ (0.01$/$0.01) & 0.0571 & 0.0219\phantom{$\mathrm{e}\mbox{-}4$} & 0.0027 &
\\[5pt]
% Mean and Cov from the truly simulated random effects:
% (1.7373,0.4988) & (0.4505,-0.0004,0.0099)
4 & $\beta$ (1) & 0.9684 & 0.1504\phantom{$\mathrm{e}\mbox{-}4$} & 0.0236 & 98\% \\
% PL: 0.2020 (51/100)
& $\sigma^2$ (0.1) & 0.0997 & 0.0023\phantom{$\mathrm{e}\mbox{-}4$} & \multicolumn{1}{c}{$5\mathrm{e}{-}6$} \\
& $\mu_1/\mu^*_1$ (2$/$1.74) & 1.7464 & 0.0460\phantom{$\mathrm{e}\mbox{-}4$} & 0.0664 & \\
& $\mu_2/\mu^*_2$ (0.50$/$0.50) & 0.4491 & 0.0948\phantom{$\mathrm{e}\mbox{-}4$} & 0.0116 & \\
& $\sigma_{11}/\sigma^*_{11}$ (0.5$/$0.45) & 0.4497 & 0.0423\phantom{$\mathrm{e}\mbox{-}4$} & 0.0043 & \\
& $\sigma_{12}/\sigma^*_{12}$ ($-$0.001$/$$-$0.0007) & -0.0439 & 0.0518\phantom{$\mathrm{e}\mbox{-}4$} &
0.0045 & \\
& $\sigma_{22}/\sigma^*_{22}$ (0.0025$/$0.0025) & 0.0797 & 0.0479\phantom{$\mathrm{e}\mbox{-}4$} &
0.0072 & \\[5pt]
% Mean and Cov from the truly simulated random effects:
% (1.7359,0.5002) & (0.4451,-0.0007,0.0025)
5 & $\beta$ (1) & 0.9934 & 0.1567\phantom{$\mathrm{e}\mbox{-}4$} & 0.0246 & 95\% \\
% PL: 0.2412 (55/100)
& $\sigma^2$ (0.1) & 0.0996 & 0.0020\phantom{$\mathrm{e}\mbox{-}4$} & \multicolumn{1}{c}{$4\mathrm{e}{-}6$}\\
& $\mu_1/\mu^*_1$ (2$/$1.74) & 1.7522 & 0.0464\phantom{$\mathrm{e}\mbox{-}4$} & 0.0636 & \\
& $\mu_2/\mu^*_2$ (0.50$/$0.50) & 0.4498 & 0.1168\phantom{$\mathrm{e}\mbox{-}4$} & 0.0162 & \\
& $\sigma_{11}/\sigma^*_{11}$ (0.50$/$0.45) & 0.4559 & 0.0442\phantom{$\mathrm{e}\mbox{-}4$} & 0.0039 &
\\
& $\sigma_{12}/\sigma^*_{12}$ ($-$0.001$/$$-$0.002) & -0.0433 & 0.0692\phantom{$\mathrm{e}\mbox{-}4$} &
0.0066 & \\
& $\sigma_{22}/\sigma^*_{22}$ (0.04$/$0.04) & 0.1186 & 0.0642\phantom{$\mathrm{e}\mbox{-}4$} & 0.0159 &
\\
\hline
% Mean and Cov from the truly simulated random effects:
% (1.7419,0.4981) & (0.4512,-0.002,0.040)
\end{tabular*}
%
% SE(MC) means the empirical standard error of
%estimates through Monte Carlo Simulation with 100 Monte Carlo samples.}
%
\end{table}

The results for the longitudinal part echo the above discussion of
the non-identifiability of the parameter $\alpha$, as the means of
the random intercept and random slopes (shown in the second column
of Table~\ref{tab1}) are consistently underestimated. The actual targets of
the estimates are the conditional quantities marked as $\mu_1^*$
and $\mu_2^*$ etc. in the first column of Table~\ref{tab1}. The sizes of
the biases vary with the level of truncation probability and size
of measurement errors and can be very small for the mean of the
random slope, for example, in setting 3, where the measurement error is
small. Thus, this bias problem in estimating the longitudinal
component may elude researchers, while it is a cause of
substantial concern in settings with large error variances.

%We consistently underestimate the mean of the random intercept, and
%lack of accuracy on the mean of the random slope.

To make statistical inference about the parameters of interest, it
is necessary to get an estimate of the standard error of the
NPMMLE, especially for $\beta$. We tried the approach in
\citet{MurpV00} and \citet{Loui82}, but neither works, so we propose
to use a bootstrap method [\citet{TsenHW05}] for estimating the
standard error of the NPMMLE and present the results in Table~\ref{tab2}.
Only the results for estimating the standard errors of
$\hat{\beta}$ and $\hat{\sigma}^2$ are shown, since they are
estimable. Table~\ref{tab2} supports the use of the bootstrap procedure, as
the estimated standard error from the bootstrap method is close to
the standard deviation from the 100 Monte Carlo samples, even when
the degree of error contamination is large or the random slopes
vary widely.
%
%[I mentioned that we used bootstrap because the method in Murphy and
%Van de Vaart and Louis method both did not work. Do we have any results
%to show? Add the references of Louis and Tseng, Hsieh and Wang (2005,
%Biometrika). ]
%
%[Jane-Ling -- The paper by Louis should be the one in 1982. I have
%changed the citation. We do have some simulation results of
%Murphy's method, including the ratio of nonnegative estimations
%out of 100 runs and the mean of those nonnegative estimations of
%SE of $\hat{\beta}$.]
%

%t2 ###
\begin{table}
\tablewidth=\textwidth
\tabcolsep=0pt
\caption{Performance of estimated variance, SE(BT), of
$\hat{\beta}$ and $\hat{\sigma}^2$ through bootstrap with 50
resamples}\label{tab2}
\begin{tabular*}{\textwidth}{@{\extracolsep{\fill}}lccc@{}}
\hline
\textbf{Case} & \textbf{Parameter} & \multicolumn{1}{c}{\textbf{SE(MC)}} & \textbf{SE(BT)}\\
\hline
1 & $\beta$ (1) & 0.1633\phantom{$\mathrm{e}\mbox{-}4$} & 0.1692\\
% PL: 0.1430 (53/100)
& $\sigma^2$ (0.1) & 0.0021\phantom{$\mathrm{e}\mbox{-}4$} & 0.0020\\[5pt]
2 & $\beta$ (1) & 0.1765\phantom{$\mathrm{e}\mbox{-}4$} & 0.1813 \\
% PL: 0.2309 (46/100)
& $\sigma^2$ (0.4) & 0.0086\phantom{$\mathrm{e}\mbox{-}4$} & 0.0091 \\[5pt]
3 & $\beta$ (1) & 0.1548\phantom{$\mathrm{e}\mbox{-}4$} & 0.1523 \\
% PL: 0.1276 (49/100)
& $\sigma^2$ (0.025) & $4.8283\mathrm{e}{-}4$ & $5\mathrm{e}{-}4$\phantom{..}\\[5pt]
4 & $\beta$ (1) & 0.1504\phantom{$\mathrm{e}\mbox{-}4$} & 0.1539 \\
% PL: 0.2020 (51/100)
& $\sigma^2$ (0.1) & 0.0023\phantom{$\mathrm{e}\mbox{-}4$} & 0.0020 \\[5pt]
5 & $\beta$ (1) & 0.1567\phantom{$\mathrm{e}\mbox{-}4$} & 0.1531\\
% PL: 0.2412 (55/100)
& $\sigma^2$ (0.1) & 0.0020\phantom{$\mathrm{e}\mbox{-}4$} & 0.0021\\
\hline
\end{tabular*}
\end{table}

%s5 ###
\section{Data example: Multi-center HIV study}\label{sec5}
In this section, we conduct an analysis on the data from a
multi-center HIV study in Italy. Details of the study design and a
previous analysis can be found in \citet{RezzLASPTSRAM89} and
\citet{Ital92}. There were 448 HIV-positive patients in the data.
The primary event of interest is the incubation period of acquired
immunodeficiency syndrome (AIDS), that is, time (in years) from
detection of HIV-infection until the onset of AIDS. There were 140
patients who received the HAART treatment at various times,
resulting in a~longitudinal treatment indicator that is fully
observable, so no modeling of this process is necessary. However,
there is a second longitudinal covariate, the CD4 counts, that are
observed only intermittently at follow-up visits, motivating the
need to model the survival and longitudinal covariates jointly.
The main biomedical interest lies in determining the effect of the
HAART treatment on reducing the risk of developing AIDS, and the
association between the incubation period of AIDS and CD4 T-cell
counts in HIV-infected subjects.

For each of the 448 subjects in the study, the longitudinal
measurements of CD4 T-cell counts were recorded intermittently
along with the time to AIDS or dropout from the study. The total
number of longitudinal measurements is 4442 and the average number
of longitudinal measurements is 9.92 per patient.

One feature of this data is that the incubation period is subject
to left-truncation and right-censoring, since patients were
recruited to the study at various times after the study began, and
only patients who have not developed AIDS at the time of
recruitment are included in the study. Moreover, only 147 out of
the 448 patients (about 33\%) developed AIDS by the end of the
study, so the right censoring rate is quite high for this data.

%
%. Since it is not
%the case that baseline measurement of CD4 counts is taken for all
%subjects (termed as delayed-entry), to make the inclusion of the
%subjects meaningful, it is required that each subject has to have
%at least one measurement of CD4 counts prior to the time at
%developing AIDS.
%
%This leads to left-truncation on the time-to-AIDS by the time at
%the first CD4 measurement.
%On the other hand, a subject may be right censored due to drop-out from
%the study or still AIDS-free by the end of the study. Thus data are
%

To model the longitudinal CD4 counts, we adopt a linear mixed
effects model on $\log(\mathrm{CD}4+1)$ with changing intercepts and
slopes at the time of HAART treatment. Thus,
\[
W^*_i(s_{ij})=X_i(s_{ij})+\varepsilon
_{ij}=A^*_{i0}+A^*_{i1}s_{ij}+A^*_{i2}I(s_{ij}>V_i)+A^*_{i3}s_{ij}I(s_{ij}>V_i)+\varepsilon
_{ij},
\]
where $\varepsilon_{ij}$ is from a normal distribution
$N(0,\sigma^2)$,
$\vec{A}^*_i=(A^*_{i0},A^*_{i1},A^*_{i2},A^*_{i3})$ is from a
4-dimensional multivariate normal distribution with a $4\times1$
mean vector $\mu$ and a $4\times4$ covariance matrix $\Sigma$, and
$V_i$ represents relative age since HIV-positive of receiving
HAART. For those who have never received HAART, $V_i$ is defined
to be infinity. For the time-to-AIDS, we assume a Cox model with
$X_i(t)$, CD4 counts, as an time-dependent covariate along with
another time-dependent treatment indicator, $I(t>V_i)$, which is
completely observed. The resulting model is
\[
\lambda(t|\vec{A}^*_i)=\lambda_0(t)\exp\bigl(\beta_1X_i(t)+\beta
_2I(t>V_i)\bigr).
\]

From the EM algorithm with Monte Carlo approximation, the slope,
$\hat{\beta}_1$, for the underlying $\log(\mathrm{CD}4 + 1)$ process
is estimated to be $-$0.5762 ($p$-value $<$ 0.001), while the slope,
$\hat{\beta}_2$, for the longitudinal treatment indicator is estimated
to be
$-$1.2189 ($p$-value $<$ 0.001). As expected, CD4 counts
are negatively associated with the risk of AIDS. One unit of
decline on $\log(\mathrm{CD}4 + 1)$ is associated with an increasing risk
of AIDS
by 78\%. In addition to its effect on CD4 counts, HAART has
an additional effect on reducing the risk of AIDS. It
significantly reduces the risk of developing AIDS by 70\% after
controlling for the CD4 counts.
%remaining parameters are as follows. The estimated variance of the
%measurement error is 0.1307. The estimated random effects
%has mean (6.5110, -0.1948, 0.1892, 0.3420) and covariance matrix
%%
%(
%%
%0.2453 & -0.0378 & 0.0183 & 0.0583 \\
%-0.0378 & 0.0498 & -0.0350 & -0.0503 \\
%0.0184 & -0.0350 & 0.1147 & 0.0058 \\
%0.0583 & -0.0503 & 0.0058 & 0.1132
%%
%).\nonumber
%%
%The HAART treatment does seem to increase CD4 counts immediately
%and has a lasting effect.\fi
Through the analysis, we confirm that
HAART effectively reduces the risk of developing AIDS both through
a positive association with patients' CD4 counts and the risk to
develop AIDS.

%s6 ###
\section{Conclusions and discussion}\label{sec6}

We have shown, both theoretically and empirically, that joint
modeling the time-to-event and longitudinal covariates is an
effective modeling approach when the time-to-event is subject to
both left truncation and right censoring. However, the extension
from right-censorship to LTRC is not trivial. By modifying the
joint likelihood, we have shown that NPMMLE leads to consistent
and asymptotically efficient estimation of the survival component
and measurement error variance under the setting of a
semiparametric Cox model. We have also demonstrated that the
corresponding EM algorithm to locate the NPMMLE has good empirical
performance and asymptotic properties under the assumption of
normal random effects. It is not only computational effective but
also robust against departures from the normality assumption.

However, one caveat is the estimability of the longitudinal
component. Although we can recover the conditional distribution of
the longitudinal parameter, $\alpha$, given $Y\geq T$, the
parameter $\alpha$ itself can not be estimated properly though the
modified likelihood due to the biased
sampling plan.
%There is not enough information available to recover the sampling bias
%induced on the longitudinal covariates.
Additional strong and possibly unverifiable assumptions might be
needed in order to recover the parameter $\alpha$ of the random
effects.
What we have accomplished in this paper is to successfully remove
the bias for the estimation of the survival components attributed
to the discrete measurement schedule and measurement errors of the
longitudinal covariates, thus permitting asymptotically valid and
efficient inference for the survival related parameters, which are
crucial for the evaluation of therapies.

\begin{appendix}

\section*{Appendix}\label{app}

%s6.1 ###
\subsection{Likelihood and the score equations}\label{appa1}
By imposing a normality assumption $N(\mu,\Sigma)$ on the random
effects $A^*_i$, the full likelihood\vadjust{\goodbreak} in (\ref{originalfulllike})
from the $i$th subject becomes
\begin{eqnarray*}
L_i^O &\propto& f_{T^*}(t_i)\sigma^{-n_i}\lambda_0(z_i)^{\delta
_i}
\\
&&{}\times\int_{-\infty}^{\infty}
\exp\Biggl\{\delta_i\beta g(z_i)a_i-\sum_{j=1}^{n_i}
[w_{ij}-g(s_{ij})a_i]^2/(2\sigma^2)\Biggr\}
Q_1(z_i,a_i)\,da_i
\\
&&{}\Big/\int_0^{\infty}\int_{-\infty}^{\infty
}Q_1(t,a_i)f_T(t)\,da_i\,dt,
\end{eqnarray*}
%
%%
%L_i^C \propto\frac{f_T(t_i)\sigma^{-n_i}\lambda_0(z_i)^{\delta_i}\int
%_{-\infty}^{\infty}
%}}[w_{ij}-g(s_{ij})a_i]^2/(2\sigma^2)\}
%Q_1(z_i,a_i)da_i}
%{\int_0^{\infty}\int_{-\infty}^{\infty
%}Q_1(t,a_i)f_T(t)da_idt},\nonumber
%%
%
where $Q_1(u,a)=\exp\{-\int_0^{u}\exp[\beta
g(t)a]\,d\Lambda_0(t)-(a-\mu)^T\Sigma^{-1}(a-\mu)/2\}$. Following
similar arguments as in \citet{Wang87} and combining with
\citet{Vard85}, we can prove that the NPMLE's of all
finite-dimensional parameters are the same as those from the
conditional likelihood of $(z_i,\delta_i,\vec{w}_i)$ given
$(Y^*_i>t_i)$. Moreover, by a proof similar to that of the
classical Cox model for right censored data, the NPMLE from the
conditional likelihood is attained by discrete baseline hazard
functions that assign positive masses only at uncensored survival
times, $(y_1^0,\ldots,y_{n_u}^0)$.

% $\lambda_0(t)$ that corresponds to a step-function, $\Lambda\{t\}$,
%whose jump sizes of a step-function $\Lambda$ at time $t$, with
%positive masses only at uncensored survival times. That is, the
%corresponding NPMLE has a discrete distribution that assigns positive
%probability mass only at uncensored survival times, $(y_1^0,

%, similar to the classical Cox model with time-independent covariates
%but with LTRC survival data.%
Let $o_i= (t_i,z_i,\delta_i, \vec{w}_i) $ denote the observed data
for the $i$th subject. The first derivative of the log full
likelihood leads to the following score functions:
\begin{eqnarray*}
s^o_{\sigma^2} &=& \sum_{i=1}^{n}
\Biggl\{\sum_{j=1}^{n_i}E[w_{ij}-A^*_ig(s_{ij})\vert
o_i]^2-n_i\sigma^2\Biggr\}\Big/\sigma^{-3},
\\
s^o_{\mu} &=& \Sigma^{-1}\sum_{i=1}^{n}E\{(A^*_i-\mu
)-[E(A^*_i\vert Y^*_i\geq T_i^*)-\mu]\vert
o_i\}
\\
&=& \Sigma^{-1}\sum_{i=1}^{n}[E(A^*_i\vert
o_i)-E(A^*_i\vert Y^*_i\geq T_i^*)],\nonumber
\\
s^o_{\Sigma} &=&
\frac{1}{2}\Sigma^{-1}\sum_{i=1}^{n}\{E[(A^*_i-\mu
)(A^*_i-\mu)^T\vert
o_i]
\\
&&\hphantom{\frac{1}{2}\Sigma^{-1}\sum_{i=1}^{n}\{}{}-E[(A^*_i-\mu)(A^*_i-\mu)^T\vert
Y^*_i\geq T^*_i]\}\Sigma^{-1},
\\
s^o_{\Lambda_k} &=& \frac{1}{\Lambda_k}-\sum_{i:t_i<y^0_k\leq
z_i}E\{\exp[\beta g(y^0_k)A^*]\vert
o_i\}-Q_2(y_k^0),
\\
s^o_{\beta} &=& \sum_{i=1}^{n}\delta
_ig(y_i)E(A^*_i\vert
o_i)
\\
&&{}-\sum_{i=1}^{n}\sum_{j:t_i<y_j^0\leq
z_i}\Lambda_jE\{g(y_j^0)A^*_i\exp[\beta g(y^0_j)A^*_i]\vert
o_i\}-Q_3,
\end{eqnarray*}
where
\begin{eqnarray*}
Q_2(y) &=& \sum_{i:y\leq t_i}E\{\exp[\beta g(y)A^*_i]\vert
o_i\}
\\
&&{}-n E\{\exp[\beta g(y)A^*_i]I(y\leq T_i)\vert Y^*_i\geq T^*_i\},
\\
Q_3 &=& \sum_{i=1}^{n}\sum_{j:y_j^0\leq
t_i}\Lambda_jE\{g(y^0_j)A^*_i\exp[\beta g(y^0_j)A^*_i]\vert
o_i\}
\\
&&{}-n E\biggl\{\sum_{j:y_j^0\leq t_i}\Lambda_jg(y_j^0)A^*_i\exp
[\beta g(y_j^0)A^*_i]\Big\vert Y^*_i\geq T^*_i\biggr\}.
\end{eqnarray*}

The score equations, $s^o_{\mu}$ and $s^o_{\Sigma}$, corresponding
to the longitudinal data reveal that the estimable terms are the
conditional mean and covariance matrix of the random effects given
that $Y^* \geq T$ rather than $\mu$ and $\Sigma$.

The score functions for $\lambda_k$,
$k=1,\ldots,n_u$, and $\beta$ have more complicated forms than
those from a partial likelihood under standard Cox model subject
to LTRC. The complication is due to the additional terms $Q_2$ and
$Q_3$, which require estimation of the expectations of nonlinear
functions of the observed data along with the the parameters of
interest. If we drop these two terms from $s^o_{\Lambda_k}$ and
$s^o_{\beta}$, the modified score functions,
$s_{\Lambda_k}=s^o_{\Lambda_k}+Q_2(y^0_k)$ and
$s_{\beta}=s^o_{\beta}+Q_3$, are exactly the score functions from
the modified likelihood. The next Lemma validates the use of the
modified likelihood (\ref{marparlike}).
% $Q_2$ and $Q_3$, and the scores $s^o_{\Lambda_k}$, $k=1,\ldots,n_u$,
%and $s^o_{\beta}$ become the counterpart from (\ref{marparlike}).

\renewcommand{\thelemma}{A.\arabic{lemma}}
\begin{lemma}\label{lem1}
%Let $s_{\Lambda_k}=s^o_{\Lambda_k}+Q_2(y^0_k)$, and $s_{\beta}=s^o_{
%
\textup{(i)}
$E_{\theta_0}(s_{\Lambda_k})=E_{\theta_0}(s^o_{\Lambda_k})$
and $E_{\theta_0}(s_{\beta})=E_{\theta_0}(s^o_{\beta})$.
This provides Fisher consistency of the estimators
(\ref{marparlike}).

\textup{(ii)} Under the regularity conditions for law of large numbers
and Slutsky theorem, $n^{-1}(s_{\Lambda_k} - s^o_{\Lambda_k})= o_p (1)$ and
$n^{-1} (s_{\beta} -s^o_{\beta})= o_p(1)$.
%o_p(1)$ asymptotically.
%
\end{lemma}

\begin{pf}
The proof follows from simple derivation and applications
of the law of large numbers along with
Slutsky's theorem.
\end{pf}

This lemma demonstrates the asymptotic equivalence of the score
functions for the survival-related parameters from
(\ref{originalfulllike}) and (\ref{marparlike}). The latter is
computationally simpler to maximize and thus more attractive than
the full likelihood.
\subsection{Proof of the consistency of the NPMMLE}
The proof of consistency includes four major steps and is
elaborated below.

\textit{Step} 1. \textit{Existence of the NPMMLE of}
($\gamma$, $\Lambda$).
We will begin the proof that the candidates for the
maximizer, $\Lambda_{n_u}$, have a finite and bounded jump at each
observed survival time. For simplicity, we use a vector form
$\vec{\lambda}_{n_u}=(\lambda_1,\ldots,\lambda_{n_u})$ to express
the jump sizes of $\Lambda_{n_u}$ at ordered survival times. The
boundedness of the jump sizes can be demonstrated by proving the
existence of an upper bound $B\in\mathbb{R}$ through apagoge.
Suppose that for any arbitrary $B\in\mathbb{R}$, there exists
$\vec{\lambda}_{n_u,B}=(\lambda_{1,B},\ldots,\lambda_{n_u,B})\in\mathbb
{R}^{n_u}\backslash[0,B]^{n_u}$
and $\gamma_B\in S_{\gamma}$ such that
$L^m(\gamma_B,\vec{\lambda}_{n_u,B})>L^m(\gamma,\vec{\lambda}_{n_u})$
for all $(\gamma,\vec{\lambda}_{n_u})\in S_{\gamma}\times
[0,B]^{n_u}$. The first part in
$L^m(\gamma_B,\vec{\lambda}_{n_u,B})$ contributed by the $i$th
subject is bounded above by
\[
(\Lambda_{n_u,B}\{z_i\}M)^{\delta_{i}}\times
\exp\biggl\{-m\sum_{j:t_{i}<y_j^0\leq
z_{i}}\lambda_{j,B}\biggr\},
\]
where $m,M$ is defined in assumption C2. Since
$\vec{\lambda}_{n_u,B}\in\mathbb{R}^{n_u}\backslash[0,B]^{n_u}$,
at least one jump size, say $\lambda_{i_0,B}$, is greater than
$B$. It induces that $\sum_{j:t_{i_0}<y_j^0\leq
z_{i_0}}\lambda_{j,B}>B$, and then implies that
$L^m(\gamma_B,\vec{\lambda}_{n_u,B})\rightarrow0$ as
$B\rightarrow\infty$. Thus $L^m(\gamma,\lambda_{n_u})=0$, for all
$(\gamma,\lambda_{n_u})\in S_{\gamma}\times\mathbb{R}^{n_u}$,
which is a contradiction. It demonstrates the boundedness of the
jump sizes of $\Lambda_{n_u}$. Along with the compactness of
$S_{\gamma}$ provided by assumption~C1, we accomplished the
existence of the NPMMLE of $(\gamma,\Lambda)$.

\textit{Step} 2. \textit{Almost surely boundedness of}
$\hat{\Lambda}(\tau_2)$ \textit{as}
$n\rightarrow\infty$.
For any fixed sample size $n$, the estimated cumulative hazard
function evaluated at the endpoint of the study can be expressed
as
%
%e6.1 ###
\begin{eqnarray}\label{consistent1}\qquad
\hat{\Lambda}(\tau_2)&=&\sum_{k=1}^{n}\frac{\delta_k
I(z_k\leq\tau_2)}{\sum_{i=1}^{n}E_{\hat{\theta}}[\exp
\{\hat{\beta}g(z_k)A^*_i\}|o_i]I(t_i<z_k\leq
z_i)}\nonumber
\\[-8pt]
\\[-8pt]
&\leq&\sum_{k=1}^{n}\frac{\delta_k
I(z_k\leq\tau_2)}{m\sum_{i=1}^{n}I(t_i<z_k\leq
z_i)}
\leq\frac{\sum_{k=1}^{n}\delta_k
I(z_k\leq\tau_2)}{m\sum_{i=1}^{n}I(t_i\leq
\tau_1)I(\tau_2\leq z_i)},\nonumber
\end{eqnarray}
where $m$ is the lower bound of $\exp\{\beta g(Y)A^*\}|Y^*\geq
T^*$, which exists under assumption C2. By the law of large numbers
and the continuous mapping theorem, we have the following two
limits as $n\rightarrow\infty$:
%
%e6.2 ###
\begin{eqnarray}
\frac{1}{n}\sum_{k=1}^{n}\delta_k
I(z_k\leq\tau_2)&\rightarrow& E\bigl(\Delta I(Y\leq\tau_2)\bigr)<1\quad
\mbox{and}\nonumber\hspace*{-35pt}
\\[-8pt]
\\[-8pt]
\frac{1}{({1}/{n})\sum_{i=1}^{n}I(t_i\leq
\tau_1\mbox{ and } \tau_2\leq
z_i)}&\rightarrow&\frac{1}{P(T\leq\tau_1\mbox{ and }Y\geq\tau_2)}<\infty
, \nonumber\hspace*{-35pt}
\end{eqnarray}
where the finiteness of the second limit is following assumption
C3. Therefore, there exists an upper bound of
$\hat{\Lambda}(\tau_2)$ even when $n$ goes to infinity. Moreover,
since the terms inside the summation in (\ref{consistent1}) are
all strictly positive, $\hat{\Lambda}(\tau_2)$ is always greater
than 0. Thus $\hat{\Lambda}(\tau_2)$ has been shown to be bounded
almost surely as $n\rightarrow\infty$.

\textit{Step} 3. \textit{Uniform convergence of} $(\hat{\sigma}^2_n,\hat
{\beta}_n,\hat{\Lambda}_n)$ \textit{to} $(\sigma^2_0,\beta_0,\Lambda_0)$.
We have shown in Step 2 that $\hat{\Lambda}(\tau_2)$ is finite,
combining with the fact that $\hat{\Lambda}$ is a~right-continuous
and nondecreasing step function along with the Helly selection
theorem, there exists a subsequence of $\hat{\Lambda}$ converging
pointwisely to a right continuous and monotone function
$\Lambda^*$ with probability 1. Moreover, by the
Balzonno--Weierstrass theorem, there is a sub-subsequence
of~$\hat{\gamma}$ which converges to some $\gamma^*$. Therefore,
there exists a sub-subsequence of~$\hat{\theta}_n$, denoted by
$\hat{\theta}_{\eta(n)}$, that converges to
$\theta^*=(\gamma^*,\Lambda^*)$. We next show that
$\theta^*=(\alpha_0^*,\sigma^2_0,\beta_0,\Lambda_0)$, where
$\alpha_0^*$ is the limit of $\hat{\alpha}$. Here a new term,
defined as
\[
\bar{\Lambda}_n(t)=\frac{1}{n}\sum_{k=1}^{n}\frac
{\delta_k
I(z_k\leq
t)}{({1}/{n})\sum_{i=1}^{n}E_{\theta_0}[\exp\{\beta
_0g(z_k)A^*_i\}|o_i]I(t_i<z_k\leq
z_i)},
\]
is introduced to serve as a
bridge between $\hat{\Lambda}_n$ and $\Lambda_0$.

We first show the convergence of $\bar{\Lambda}_n$ to $\Lambda_0$
as follows. We will use a~property that the class of all functions
from a closed set to $\mathbb{R}$, which are uniformly bounded and
of bounded variation, is Glivenko--Cantelli. Consider the
denominator of $\bar{\Lambda}_n$. The assumptions imply that
functions of the form $u\rightarrow
E_{\theta_0}[\exp\{\beta_0g(u)A^*\}I(T<u\leq Y)|o]$, where $o$
denotes the observed data of a subject, are uniformly bounded and
of bounded variation, so the class of these functions is
Glivenko--Cantelli. Therefore,
%
%e6.3 ###
\begin{eqnarray}\label{consistent2}
&&\frac{1}{n}\sum_{i=1}^{n}E_{\theta_0}[\exp\{\beta
_0g(u)A^*_i\}|o_i]I(t_i<u\leq
z_i)\nonumber
\\[-8pt]
\\[-8pt]
&&\qquad\rightarrow E_{\theta_0}\{\exp[\beta_0g(u)A^*]I(T<u\leq
Y)|Y^*\geq T^*\}\nonumber
\end{eqnarray}
uniformly on $[\tau_1,\tau_2]$. Along with assumption C4, the
uniform convergence of the inverse of the right-hand side to the
inverse of the left-hand side in (\ref{consistent2}) holds.
Moreover, uniform boundedness and bounded variation of functions
$t\rightarrow\Delta I(Y\leq t)$ imply the Glivenko--Cantelli
property of the class consisting of them. Thus we also have
%
%e6.4 ###
\begin{equation}\label{consistent3}
\frac{1}{n}\sum_{i=1}^{n}\Delta_iI(Y_i<t)\rightarrow
E_{\theta_0}[\Delta I(Y<t)]
\end{equation}
uniformly on $[\tau_1,\tau_2]$. Since
\[
\Lambda_0(t)=E\biggl[\frac{\Delta
I(Y\leq t)}{E\{\exp[\beta_0g(Y)A^*]I(T<u\leq Y)|Y^*\geq
T^*\}|{u=Y}}\biggr],
\]
combining the convergence of the inverse of both
sides in (\ref{consistent2}) and (\ref{consistent3}), we obtain
$\bar{\Lambda}_n$ converges uniformly to $\Lambda_0$ on
$[\tau_1,\tau_2]$. By considering the uniform convergence of the
ratio of $\hat{\Lambda}\{u\}/\bar{\Lambda}\{u\}$ to
$d\Lambda^*(u)/d\Lambda_0(u)$ for $u \in[\tau_1,\tau_2]$, as
demonstrated on pages 2146--2147 in \citet{ZengC05}, the uniform
convergence of $\hat{\Lambda}$ to $\Lambda^*$ is established. The
remaining task is to prove the equivalence of
$\theta^*=(\gamma^*,\Lambda^*)$ and
$\theta_0^*=(\beta_0,\sigma^0,\alpha^*,\Lambda_0)$. This can be
done by considering the empirical mean of the distance between
$l_i^{m}(\hat{\theta}_n)$ and
$l_i^{m}(\beta_0,\sigma^0,\alpha^*,\bar{\Lambda}_n)$ and
demonstrating that
$\mathrm{E}_{\theta_0^*}[l^{m}(\theta^*)/l^{m}(\theta_0^*)]=0$
almost surely as shown on page 910 in \citet{DupuGM06}. Thus
$(\hat{\sigma}^2,\hat{\beta},\hat{\Lambda_0})$ converges uniformly
to $(\sigma^2_0,\beta_0,\Lambda_0)$.

%s6.3 ###
\subsection{Proof of asymptotic normality of the NPMMLE}
We will apply Theorem 3.3.1 in \citet{VaarW96} to prove the
asymptotic normality of the NPMMLE ($\hat{\gamma}$,
$\hat{\Lambda}$). The proof consists of four steps to verify each
of the four conditions in
their theorem.

\textit{Step} 1. \textit{Fr\'{e}chet differentiability of the score functions}.
%
%The modified likelihood contributed by the $i$th subject can be
%re-written as
%%
%L_i^{m}(\gamma,\Lambda) = \int&[\Lambda\{y_i\}]^{\delta_i}\exp\{\delta
%_i\beta a_i^Tg(y_i)-\overset{n_u}{\underset{k=1}{\sum}}\Lambda\{y_k^0\}
%&(2\pi
%}}[w_{ij}-a_i^Tg(t_{ij})]^2/(2\sigma^2)\}
%f_A(a_i) da_i,\nonumber
%
%
For notation simplification, the parameter $\sigma^2$
will be combined with $\alpha$ into $\gamma_1= (\sigma^2, \alpha)$
so that the single parameter $\gamma_1$ denotes the parameter of
the measurement error $\varepsilon$ and the latent random variable
$A^*$. Thus, the new parameter vector is
$\theta=(\gamma_1,\beta,\Lambda)$.

Consider a one-dimensional submodel along the direction
$(h_1,h_2,h_3)$ of the form
\[
\theta_t=\bigl(\gamma_1+th_1,\beta+th_2,\Lambda_t(h_3)\bigr),
\]
where
\[
\Lambda_t(h_3)(\cdot)=\int_0^{\cdot}\bigl(1+th_3(u)\bigr)\,d\Lambda(u),
\]
$h_1\in\mathbb{R}^d$, $h_2\in\mathbb{R}$ and $h_3$ is a
bounded-variation function on $[0,\tau_2]$. Let
$H=\{h=(h_1,h_2,h_3)\}$ and $H_p=\{h\in H\dvtx\Vert h_1\Vert,
\vert h_2\vert, \Vert h_3\Vert_v\leq p\}$. The notation
\mbox{$\Vert\cdot\Vert_v$} denotes the absolute value evaluated
at 0 plus the total variation of the argument. The imputed
log-likelihood contributed by the $i$th subject evaluated at
$\theta$, given the current value of parameter denoted as
$\tilde{\theta}$, is denoted by~$l_{\tilde{\theta},i}(\theta)$.
The corresponding score function of the local parameter $t$ is
\begin{eqnarray*}
\frac{\partial}{\partial t}l_{\tilde{\theta},i}(\theta_t)
&=&h_2\mathrm{E}_{\tilde{\theta}}\biggl[\delta_ig(z_i)A^*_i
\\
&&\hphantom{h_2\mathrm{E}_{\tilde{\theta}}\biggl[}{}-\int
_{t_i}^{z_i}g(u)A^*_i\exp\{(\beta+th_2)g(u)A^*_i\}\bigl(1+th_3(u)\bigr)\,d\Lambda
(u)\Big|o_i\biggr]
\\
&&{}-\mathrm{E}_{\tilde{\theta}}\biggl[\int_{t_i}^{z_i}h_3(u)\exp\{(\beta
+th_2)g(u)A^*_i\}\,d\Lambda(u)\Big|o_i\biggr]
\\
&&{}+h_1\mathrm{E}_{\tilde{\theta}}\biggl[\frac{\partial}{\partial
t}f_{\varepsilon,A^*}(\varepsilon,A^*_i|\gamma_1+th_1)\Big|o_i\biggr]
\\
&&{}+\frac{\delta_ih_3(z_i)}{1+th_3(z_i)}.
\end{eqnarray*}
Thus the imputed score function of $t$ contributed by the $n$
subjects evaluated at $t=0$ is
%
%e6.5 ###
\begin{eqnarray}\label{normal1}
S_{n,\tilde{\theta}}(\theta)(h)&=&\frac{1}{n}\sum_{i=1}^{n}\frac{\partial}{\partial
t}l_{\tilde{\theta},i}(\theta_t)\bigg\vert_{t=0}\nonumber
\\[-9pt]
\\[-9pt]
&=&h_1^T S_{n,\tilde{\theta},1}(\theta)+h_2
S_{n,\tilde{\theta},2}(\theta)+S_{n,\tilde{\theta},3}(\theta
)(h_3),\nonumber
\end{eqnarray}
where
\begin{eqnarray*}
S_{n,\tilde{\theta},1}(\theta)&=&\frac{1}{n}\sum_{i=1}^{n}
\mathrm{E}_{\tilde{\theta}}\biggl[\frac{\partial}{\partial\gamma_1}\log
f_{\varepsilon,A^*}(\varepsilon_i,A^*_i|\gamma_1)\Big|o_i\biggr],
\\[-1pt]
S_{n,\tilde{\theta},2}(\theta)&=&\frac{1}{n}\sum_{i=1}^{n}
\mathrm{E}_{\tilde{\theta}}\biggl[\delta_ig(z_i)A^*_i-\int
_{t_i}^{z_i}g(u)A^*_i\exp\{\beta
g(u)A^*_i\}\,d\Lambda{u}\Big|o_i\biggr],
\\[-1pt]
S_{n,\tilde{\theta},3}(\theta)(h_3)&=&\frac{1}{n}\sum_{i=1}^{n}
\biggl\{\delta_i
h_3(z_i)-\mathrm{E}_{\tilde{\theta}}\biggl[\int_{t_i}^{z_i}h_3(u)\exp\{\beta
g(u)A^*_i\}\,d\Lambda(u)\Big|o_i\biggr]\biggr\}.
\end{eqnarray*}
By defining
$\theta(h)\,{=}\,(\gamma_1,\beta,\Lambda)(h_1,h_2,h_3)\,{=}\,h_1^T\gamma_1\,{+}\,h_2\beta
\,{+}\,\int_{0}^{\tau2}h_3(u)\,d\Lambda(u)$,
where $h\in H_p$, the parameter $\theta$ can be regarded as a
functional on $H_p$, the parameter space $\Theta=\{\theta\}$ is a
subspace of $L^{\infty}(H_p)$ and the score in (\ref{normal1}) is
a random map from $\Theta$ to a Banach space which contains
functions (operations) of $h$.

Besides the above imputed score, we also need the mean imputed
score function of $t$ under the true value $\theta_0$ and denote
it as
\[
S_{\tilde{\theta}}(\theta)(h)=\mathrm{E}_{\theta_0}\biggl[\frac{\partial
}{\partial
t}l_{\tilde{\theta}}(\theta_t)\bigg\vert_{t=0}\biggr]=h_1^T
S_{n\tilde{\theta},1}(\theta)+h_2
S_{\tilde{\theta},2}(\theta)+S_{\tilde{\theta},3}(\theta)(h_3),
\]
where
\begin{eqnarray*}
S_{\tilde{\theta},1}(\theta)
&=&\mathrm{E}_{\theta_0}\biggl\{
\mathrm{E}_{\tilde{\theta}}\biggl[\frac{\partial}{\partial\gamma_1}\log
f_{\varepsilon,A^*}(\varepsilon_i,A^*_i|\gamma_1)\Big|o_i\biggr]\biggr\},
\\[-1pt]
S_{\tilde{\theta},2}(\theta)
&=&\mathrm{E}_{\theta_0}\biggl\{
\mathrm{E}_{\tilde{\theta}}\biggl[\Delta_ig(Y_i)A^*_i-\int
_{T_i}^{Y_i}g(u)A^*_i\exp\{\beta
g(u)A^*_i\}\,d\Lambda(u)\Big|o_i\biggr]\biggr\},
\\[-1pt]
S_{\tilde{\theta},3}(\theta)(h_3)
&=&\mathrm{E}_{\theta_0} \biggl\{\Delta_i
h_3(Y_i)-\mathrm{E}_{\tilde{\theta}}\biggl[\int_{T_i}^{Y_i}h_3(u)\exp\{\beta
g(u)A^*_i\}\,d\Lambda(u)\Big|o_i\biggr]\biggr\}.
\end{eqnarray*}

To prove the Fr\'{e}chet differentiability of the map, $\theta
\rightarrow S_{\theta^*_0}(\theta)$ at $\theta^*_0$, where
$\theta^*_0=(\gamma^*_{10},\beta_0,\Lambda_0)$ with
$\gamma^*_{10}=(\sigma_0^2,\alpha^*)$, we need to calculate the
corresponding derivative. First, we introduce a notation
$\nabla_{\theta}S_{\tilde{\theta}}(\theta^*_0)=\frac{\partial}{\partial
t}S_{\tilde{\theta}}(\theta^*_0+t\theta)\vert_{t=0}$, where
$\theta^*_0+t\theta=(\alpha_0+t\alpha,\beta_0+t\beta,\Lambda_0(\cdot
)+t\Lambda(\cdot))$.
Then
%
%e6.6 ###
\begin{eqnarray}\label{normal2}
&&\nabla_{\theta}S_{\tilde{\theta}}(\theta^*_0)(h)\nonumber
\\[-1pt]
&&\quad=\frac{\partial
}{\partial
t}S_{\tilde{\theta}}(\theta^*_0+t\theta)(h)\bigg\vert_{t=0}\nonumber
\\[-1pt]
&&\quad=\frac{\partial}{\partial t}
\mathrm{E}_{\theta^*_0}\biggl\{h_1^T\mathrm{E}_{\tilde{\theta}}\biggl[\frac{\partial
}{\partial(\gamma^*_{10}+t\gamma_1)}\log
f_{\varepsilon,A^*}(\varepsilon_i,A^*_i|\gamma^*_{10}+t\gamma_1)\Big|o_i\biggr]\nonumber
\\[-1pt]
&&\hphantom{\quad=\frac{\partial}{\partial t}
\mathrm{E}_{\theta^*_0}\biggl\{}{}+h_2\mathrm{E}_{\tilde{\theta}}\biggl[\Delta_ig(Y_i)A^*_i\nonumber
\\[-9pt]
\\[-9pt]
&&\hphantom{\quad=\frac{\partial}{\partial t}
\mathrm{E}_{\theta^*_0}\biggl\{+h_2\mathrm{E}_{\tilde{\theta}}\biggl[}{}-\int\!
_{T_i}^{Y_i}g(u)A^*_i\exp\{\beta_0+t\beta
g(u)A^*_i\}\bigl(d\Lambda_0(u)+t\,d\Lambda(u)\bigr)\Big|o_i\biggr]\nonumber
\\[-1pt]
&&\hphantom{\quad=\frac{\partial}{\partial t}
\mathrm{E}_{\theta^*_0}\biggl\{}{}+\Delta_ih_3(Y_i)\nonumber
\\[-1pt]
&&\hphantom{\quad=\frac{\partial}{\partial t}
\mathrm{E}_{\theta^*_0}\biggl\{}{}
-\mathrm{E}_{\tilde{\theta}}\biggl[\int
_{T_i}^{Y_i}h_3(u)\nonumber
\\[-1pt]
&&\hphantom{\quad=\frac{\partial}{\partial t}
\mathrm{E}_{\theta^*_0}\biggl\{-\mathrm{E}_{\tilde{\theta}}\biggl[\int}{}\times\exp\{(\beta_0+t\beta)
g(u)A^*_i\}\bigl(d\Lambda_0(u)+t\,d\Lambda(u)\bigr)\Big|o_i\biggr]\biggr\}\bigg\vert_{t=0}.\nonumber
\end{eqnarray}
Using the chain rule, equation (\ref{normal2}) can be simplified as
\[
-\gamma_1^T\sigma_{\tilde{\theta},1}(h)-\beta\sigma_{\tilde{\theta
},2}(h)-\int_{0}^{\tau_2}\sigma_{\tilde{\theta},3}(h)(u)\,d\Lambda
(u),\nonumber
\]
where
%
%e6.9 ###
%e6.8 ###
%e6.7 ###
\begin{eqnarray}
\label{normal3}
\hspace*{40pt}\sigma_{\tilde{\theta},1}(h)
&=&-\mathrm{E}_{\theta^*_0}\biggl\{h_1^T\mathrm
{E}_{\tilde{\theta}}\biggl[\frac{\partial^2}{\partial\gamma_1\,\partial\gamma_1^T}
\log f_{\varepsilon,A^*}(\varepsilon_i,A^*_i\vert
\gamma^*_{10})\Big\vert o_i\biggr]\biggr\},
\\
\label{normal4}\sigma_{\tilde{\theta},2}(h)
&=&\mathrm{E}_{\theta^*_0}\biggl\{\mathrm
{E}_{\tilde{\theta}}\biggl[\int_{0}^{\tau_2}
[h_2g(u)A^*_i+h_3(u)]g(u)A^*_i\exp\{\beta_0g(u)A^*_i\}\nonumber
\\[-8pt]
\\[-8pt]
&&\hspace*{121pt}{}\times I(T_i<u\leq Y_i)\,d\Lambda_0(u)\Big\vert o_i\biggr]\biggr\},\nonumber
\\
\label{normal5}\sigma_{\tilde{\theta},3}(h)(u)&=&\mathrm{E}_{\theta^*_0}\bigl\{\mathrm
{E}_{\tilde{\theta}}\bigl[
[h_2g(u)A^*_i+h_3(u)]\exp\{\beta_0g(u)A^*_i\}\nonumber
\\[-8pt]
\\[-8pt]
&&\hspace*{105pt}{}\times I(T_i<u\leq Y_i)\vert o_i\bigr]\bigr\}.\nonumber
\end{eqnarray}
Evaluating (\ref{normal2}) at the true value $\theta^*_0$ leads to
%
%e6.10 ###
\begin{equation}\label{normal6}
\nabla_{\theta}S_{\theta^*_0}(\theta^*_0)(h)=
-\gamma_1^T\sigma_{\theta^*_0,1}(h)-\beta\sigma_{\theta^*_0,2}(h)-\int
_{0}^{\tau_2}\sigma_{\theta^*_0,3}(h)(u)\,d\Lambda(u),\hspace*{-35pt}
\end{equation}
where each of the $\sigma$-function has similar forms as the
corresponding function listed in
(\ref{normal3}), (\ref{normal4}) or (\ref{normal5}) with the double
expectation
$\mathrm{E}_{\theta^*_0}\{\mathrm{E}_{\tilde{\theta}}\{[\cdot\vert
o_i]\}$ replaced by $\mathrm{E}_{\theta^*_0}\{\cdot\}$. Now apply
the Taylor expansion of $\exp\{(\beta_0+t\beta)g(u)A^*_i\}$ at
$t=0$, to get
\[
S_{\theta^*_0}(\theta^*_0+t\theta)-S_{\theta^*_0}(\theta^*_0)-\nabla
_{\theta}S_{\theta^*_0}(\theta^*_0)=o(t),
\]
where the small-$o$ function does not depend on $\theta$. Therefore,
\[
\frac{\Vert
S_{\theta^*_0}(\theta^*_0+t\theta)-S_{\theta^*_0}(\theta^*_0)-\nabla
_{\theta}S_{\theta^*_0}(\theta^*_0)\Vert_{p}}{t}
\rightarrow0\qquad\mbox{as }t\rightarrow0
\]
uniformly in $\theta=(\gamma_1,\beta,\Lambda)$. Thus the
Fr\'{e}chet derivative of the mapping $\theta\rightarrow
S_{\theta^*_0}(\theta)$ evaluated at $\theta^*_0$ takes the form
(\ref{normal6}). We will use the notation
$\dot{S}_{\theta^*_0}(\theta^*_0)(\theta)$ to denote it.

\textit{Step} 2. \textit{Continuous invertibility of} $\dot{S}_{\theta
^*_0}(\theta^*_0)(\theta)$.
The continuous invertibility of the Fr\'{e}chet derivative can be
established by showing that there exists some number $c>0$ such
that
%
%e6.11 ###
\begin{equation}\label{normal7}
\inf_{\theta\in\mathrm{lin}\Theta}\frac{\Vert\dot
{S}_{\theta^*_0}(\theta^*_0)\Vert_{l^\infty(H)}}{\Vert\theta
\Vert_{l^\infty(H)}}>c.
\end{equation}
Since $\dot{S}_{\theta^*_0}(\theta^*_0)(\theta)$ can be expressed
as a linear combination of the three $\sigma$-operators according
to (\ref{normal2}), it is necessary to check the continuous
invertibility of those $\sigma$-operators. The proof is similar to
the arguments in the Appendix of \citet{ZengC05}. Through the
continuous invertibility of~$\sigma_{\theta^*_0}$, the lower bound
$c$ can be found as $\frac{q}{3p}$, where $q$ satisfies
$\sigma_{\theta^*_0}^{-1}(H_q)\subseteq H_p$. Details to find the
lower bound are analogous to the approach in \citet{DupuGM06}
(page 915). Thus the derivative $\dot{S}_{\theta^*_0}(\theta^*_0)$
is continuously invertible.

\textit{Step} 3. \textit{Convergence in distribution to a tight element}.
In this step, the convergence of
$\sqrt{n}(S_{n,\hat{\theta}_n}-S_{\theta^*_0})(\theta^*_0)$ in
distribution will be demonstrated. Since
$S_{\theta^*_0}(\theta^*_0)$ is the mean of the score function
evaluated at the true value of $\theta$, it is equal to zero. Then
\[
[S_{n,\hat{\theta}_n}-S_{\theta^*_0}](\theta^*_0)(h)=\frac{1}{n}\sum_{i=1}^{n}[
D_{i,1}(h)+D_{i,2}(h)+\delta_i h_3(y_i)+D_{i,3}(h)],
\]
where
\begin{eqnarray*}
D_{i,1}(h)
&=&h_1^T\mathrm{E}_{\hat{\theta}_n}\biggl[\frac{\partial}{\partial
\alpha}\log
f_{\varepsilon,A^*}(\varepsilon_i,A^*_i|\gamma^*_{10})\Big|o_i\biggr],
\\
D_{i,2}(h)
&=&h_2\mathrm{E}_{\hat{\theta}_n}\biggl[\delta_ig(z_i)A^*_i-\int
_{t_i}^{z_i}g(u)A^*_i\exp\{\beta_0
g(u)A^*_i\}\,d\Lambda_0{u}\Big|o_i\biggr],
\\
D_{i,3}(h)
&=&-\mathrm{E}_{\hat{\theta}_n}\biggl[\int_{t_i}^{z_i}h_3(u)\exp\{
\beta_0
g(u)A^*_i\}\,d\Lambda_0(u)\Big|o_i\biggr].
\end{eqnarray*}
The class $\{\frac{1}{n}\sum(D_{i,1}+D_{i,2})(h)\dvtx\Vert
h_1\Vert+\vert h_2\vert\leq p\}$ is bounded Donsker, since it is
a finite dimensional class of measurable score functions.
Moreover, since any class of real-valued functions on $[0,\tau_2]$
that are uniformly bounded and bounded in variation is Donsker,
the class $\{\delta h_3(y)\dvtx h_3\in BV_p\}$ is Donsker. The Donsker
property of the class $\{\frac{1}{n}\sum D_{i,3}(h)\dvtx h_3\in
BV_p\}$ also follows from this fact. We have thus shown that the
class
$\{[S_{n,\hat{\theta}_n}-S_{\theta^*_0}](\theta^*_0)(h)\dvtx\Vert
h_1\Vert+\vert h_2\vert\leq p,h_3\in BV_p\}$ is Donsker, since
the sum of bounded Donsker classes is also Donsker. This implies
\[
\sqrt{n}(S_{n,\hat{\theta}_n}-S_{\theta^*_0})(\theta^*_0)\mathop
{\rightarrow}^{D}Z,
\]
a tight Gaussian process in $l_{\infty}(H_p)$.

\textit{Step} 4. \textit{Verification of conditions} 1 \textit{and} 4.
Condition 4 holds by the consistency of the estimator
$\hat{\theta}_n$. Condition 1 can be verified by considering the
Donsker property of the class
$\{S_{\cdot,\theta}(\theta)(h)-S_{\cdot,\theta^*_0}(\theta
^*_0)(h)\dvtx\Vert
\theta-\theta^*_0\Vert_p<\nu,h\in H_p\}$ for some $\nu>0$,
where $S_{\cdot,\theta}(\theta)(h)$ is the general form of
$S_{i,\theta}(\theta)(h)=\frac{\partial}{\partial
t}l_{\theta,i}(\theta_t)\vert_{t=0}$. We omit the details since
they are similar to those for the case of right-censored data,
considered in \citet{ZengC05}.

We have verified the four conditions needed for the asymptotic
distribution of the NPMMLE $\hat{\theta}_n$, and therefore
\[
\sqrt{n}(\hat{\theta}_n-\theta^*_0)\mathop{\rightarrow}^{D}-\dot
{S}_{\theta^*_0}(\theta^*_0)Z,
\]
as $n\rightarrow\infty$.

Using the form of the Fr\'{e}chet derivative in (\ref{normal2}),
one finds that there exists a linear operator
$\sigma=(\sigma_{\theta^*_0,1},\sigma_{\theta^*_0,2},\sigma_{\theta^*_0,3})$
that maps $H_p$ to $\mathbb{R}^{d+1}\times BV_p$, such that
\begin{eqnarray*}
\dot{S}_{\theta^*_0}(\theta^*_0)(\theta_1-\theta_2)(h)
&=&
-(\gamma_{11}-\gamma_{12})^T\sigma_{\theta^*_0,1}(h)-(\beta_1-\beta
_2)\sigma_{\theta^*_0,2}(h)
\\
&&{}-\int_{0}^{\tau2}\sigma_{\theta^*_0,3}(h)(u)\,d(\Lambda_1-\Lambda
_2)(u).\nonumber
\end{eqnarray*}
The continuous invertibility of the $\sigma$ operator has been
shown already, so its inverse operator, denoted by $\sigma^{-1}$,
exists. Since
\[
\sqrt{n}\dot{S}_{\theta^*_0}(\theta^*_0)(\hat{\gamma}_1-\gamma
^*_{10},\hat{\beta}-\beta_0,\hat{\Lambda}-\Lambda_0)(h)
=\sqrt{n}\{S_{n,\theta^*_0}(h)-S_{\theta^*_0}(\theta^*_0)(h)\}
+o_p(1),
\]
by applying the inverse operator $\sigma^{-1}$ on both sides we
obtain that
%
%e6.12 ###
\begin{eqnarray}\label{normal8}
&&\sqrt{n}\biggl\{-(\hat{\gamma}_1-\gamma^*_{10})^Th_1-(\hat{\beta}-\beta
_0)h_2-\int_{0}^{\tau2}h_3(u)\,d(\hat{\Lambda}-\Lambda_0)(u)\biggr\}\nonumber
\\[-8pt]
\\[-8pt]
&&\qquad=\sqrt{n}\{S_{n,\theta^*_0}(\tilde{h})-S_{\theta^*_0}(\theta
^*_0)(\tilde{h})\}+o_p(1),\nonumber
\end{eqnarray}
where
$\tilde{h}=(\tilde{h}_1,\tilde{h}_2,\tilde{h}_3)=\sigma^{-1}(h)$.
If $h_1$ and $h_3$ in (\ref{normal8}) are chosen to be 0, then this
reduces to
\begin{eqnarray*}
&&\sqrt{n}\{-(\hat{\beta}-\beta_0)h_2\}
\\
&&\qquad=\sqrt{n}\{S_{n,\theta^*_0}(\tilde{h})-S_{\theta^*_0}(\theta
^*_0)(\tilde{h})\}+o_p(1),
\end{eqnarray*}
where the latter term is in the form of linear combinations of
score functions for the parameters. Since score functions derived
from the modified likelihood is asymptotically equivalent to those
from the full likelihood by Lemma~\ref{lem1}, the influence function is the
same as the efficient influence function for $\beta_0h_2$ by its
uniqueness in the linear span of the scores. Thus
the estimator $\hat{\beta}$ is efficient for $\beta_0$.
\end{appendix}

\section*{Acknowledgments}
The authors thank the Associate Editor and reviewers for
insightful comments.

% imsref loaded by audrone.aklyte, 2012-05-15 14:47:44
% imsref loaded by audrone.aklyte, 2012-05-15 15:41:15

%suskaldyti doi

\printaddresses

\end{document}